\documentclass[11pt, letter]{article}
\usepackage[letterpaper, left=1.truein, right=.8truein, top = .5truein, bottom = 1.5truein]{geometry}

\usepackage{lineno}

\usepackage{eulervm}
\usepackage{tgpagella}
\usepackage[T1]{fontenc}
\usepackage{amsmath}
\usepackage{amssymb}
\usepackage{amsthm}
\theoremstyle{definition}
\usepackage{cite}
\usepackage{color}
\usepackage{graphicx}
\usepackage{caption}
\usepackage{color}
\usepackage{bm}

\usepackage{scalerel}

\DeclareMathAlphabet{\mathpzc}{OT1}{pzc}{m}{it}

\usepackage{algorithm}
\usepackage{algpseudocode}
\makeatletter \def\BState{\State\hskip-\ALG@thistlm} \makeatother

\usepackage{etoolbox} \makeatletter
\patchcmd{\@maketitle}{\begin{center}}{\begin{flushleft}}{}{}
\patchcmd{\@maketitle}{\begin{tabular}[t]{c}}{\begin{tabular}[t]{@{}l}}{}{}
\patchcmd{\@maketitle}{\end{center}}{\end{flushleft}}{}{}

\newcommand\be{\begin{equation}}
\newcommand\ee{\end{equation}}

\newcommand{\zero}{{\bm{0}}}

\newtheorem{thm}{Theorem}[section]
\newtheorem{lem}[thm]{Lemma}
\newtheorem{prp}[thm]{Proposition}
\newtheorem{cor}[thm]{Corollary}
\newtheorem{dfn}[thm]{Definition}
\newtheorem{rem}[thm]{Remark}
\newtheorem{obs}[thm]{Observation}

\newcommand{\Real}{\mathbb{R}}

\newcommand{\tU}{\tilde{U}}
\newcommand{\tf}{\tilde{f}}

\def\Real{\mathbb{R}}

\def\ex{\mathbb{E}}

\def\prob{\mathbb{P}}
\def\quadr{\mathtt{Q}}

\def\mtree{\mathtt T}

\def\fra{\mbox{\textsf{{\textbf{N}}}}}
\def\maa{\reflectbox{\fr}}

\def\fr{\fra\ }
\def\ma{\maa\ }

\newcommand\Em{\mathbf{e}}
\newcommand\Eb{\mathbf{b}}
\newcommand\Ed{\mathbf{d}}

\newcommand\cov{\mathtt{Cov}}
\newcommand\wi{\mathbf{w}}

\newcommand{\sif}{\mathcal{F}}
\newcommand{\ph}{\mathbf{PH}}
\newcommand{\brm}{\mathbf{B}}
\newcommand{\brmd}[1]{\brm^{#1}}
\newcommand{\mph}{\bm{\mu}}
\newcommand\path{\bm{\pi}}
\newcommand{\arpath}[3]{{#1}\xrightarrow{\displaystyle{\scaleobj{.5}{#2}}}{#3}}
\newcommand\vertical{\mathtt{B}}

\def\Push{\mathtt{Push}}
\def\Pop{\mathtt{Pop}}
\def\Top{\mathtt{Top}}

\def\Bot{\mathtt{Bot}}
\def\Dir{\mathtt{Dir}}

\usepackage[utf8]{inputenc}
\usepackage[russian,english]{babel}

\newcommand{\inl}{{\bm{\alpha}}}
\newcommand{\intervals}{\bm{A}}
\newcommand{\boundaries}{\bm{BA}}

\def\Push{\mathtt{Push}}
\def\Pop{\mathtt{Pop}}
\def\Top{\mathtt{Maxima}}

\def\Bot{\mathtt{Minima}}
\def\Dir{\mathtt{Direction}}

\def\pth{\bm{\pi}}

\usepackage{algorithm}
\usepackage{algpseudocode}

\makeatletter \def\BState{\State\hskip-\ALG@thistlm} \makeatother

\usepackage{etoolbox} \makeatletter
\patchcmd{\@maketitle}{\begin{center}}{\begin{flushleft}}{}{}
\patchcmd{\@maketitle}{\begin{tabular}[t]{c}}{\begin{tabular}[t]{@{}l}}{}{}
\patchcmd{\@maketitle}{\end{center}}{\end{flushleft}}{}{}

\usepackage{mathtools}
\AtBeginDocument{}

\begin{document}


\title{Time Series, Persistent Homology and Chirality} 

\date{\today}

\author{Yuliy Baryshnikov\footnote{Partially supported by NSF via
  grant DMS-1622370 and via MURI {\em SLICE}.}\\ Departments of Mathematics and ECE,\\ University
  of Illinois at Urbana-Champaign 
  \\ 
  {\tt
    ymb@illinois.edu} }

\maketitle

\abstract{Interactions of the maxima and minima of the univariate functions appear in combinatorics as Dyck paths, and in topological data analysis as persistent homology. We study these descriptors for Brownian motions with drift, deriving the intensity measure and correlation functions for the persistence diagram point process $\ph_0$, and quantifying the intrinsic asymmetries in the coupling of maxima and minima.}

\section{Introduction}
Motivation for this study comes from a growing body of work addressing Morse theory of random function, in the context of {\em Topological Data Analysis} (see, e.g. \cite{adler_taylor, random_comp}).

Morse theory associates to a smooth function on a manifold a
collection of its critical points with some of their local descriptors (such
as critical values or indices, for nondegenerate
critical points). It became a powerful tool to characterize
both the function and the underlying manifold in different contexts,
such as spectral theory of differential operators, or diffusions on
manifolds.

An application of the Morse descriptors in {data analysis} emerged
over the past decade, popularized under the name of {\em persistent
  homology}. The original motivation behind the notion was to provide
a scale-independent toolbox for understanding the topology of the unknown
underlying model for the data; the
properties of the function defining the persistence data was a mere byproduct. With time, however, persistent homology became
a powerful tool to sketch the properties of the function itself.

The analysis of the mapping that associates to a function its persistent homology is rather subtle, and it is natural to attempt to understand its output on some "typical functions" proxied, customarily, by realizations of a {\em random function}. While in general precise characterization of the output of the persistent homology black box is rather implicit (compare, for example, \cite{adler_topological_2011,kahle_limit_2010, bobrowski_distance_2011}), there is a class of random functions where that can be described quite precisely: namely, the trajectories of certain Brownian motions. This is the central object of study of this paper. We use the standard properties of Brownian motions to investigate in details 0-dimensional persistent homology, interpreted as a persistence diagram point process. In particular, we find the intensity density and 2-point correlation functions for the standard Brownian motions with drift.

While understanding the structure of persistent homology for random univariate functions was the starting point of this paper, our findings address a somewhat finer descriptor than just persistence diagram. As we discuss below, for univariate functions the zeroth persistent homology is just a sketch of a finer invariant of the function, its {\em merge tree}. In particular, each bar comes with a right- or left- (\fra- or \maa-, in our nomenclature) orientation, a {\em chirality}. We find relative frequencies of them. 


Main qualitative results of this paper are the formulae for the density of the $\ph_0$ point process and its 2-point correlation functions (Proposition \ref{prp:intensity} and Theorem \ref{thm:g2}) and the results about the expected excess of the number of \fr bars over the number of \ma bars, Theorem \ref{thm:foverm} for Brownian motions with positive drift.

The results of this paper should be seen in the general context of {\em reparametrization invariant} tools of data analysis, in our situation of analysis of time series. Changing the coordinates of the  underlying space corrupts many traditional tools (such as Fourier analysis). By the contrast, persistence diagrams, their various refinements described in this paper, or the characteristics like unimodal category \cite{baryshnikov_minimal_2018} survive right compositions with homeomorphism, and therefore describe patterns independent of the semantics attached to the domain coordinates.

The structure of the paper: in Section \ref{sec:ph} we recall basics of persistent homology, merge trees; specialize these definitions to the univatiate case, and introduce our notion of chirality. In Section \ref{sec:brownian} we introduce the random persistence diagrams for trajectories of Brownian motions, and use that in Section \ref{sec:chir_brown} to evaluate the average numbers of bars of different chiralities in the Brownian motions with constant drift. In Section \ref{sec:covar} we develop the formalism of automata, and use it to compute the 2-point correlation functions for $\ph_0$.

In Appendix we present, inter alia, an efficient algorithm for computing $\ph_0$ for time series, and introduce a monoid structure on merge trees.

\section{Persistent Homology and Merge Trees}\label{sec:ph}
\subsection{Persistent Homology Diagrams}
Here we briefly mention some basic notions of the persistent homology, to place our results in the appropriate context.

Consider a continuous function $f:X\to\Real$ on a topological space $X$. (In most applications one can assume that $X$ can be triangulated.)
Sublevel sets of $f$ define a filtration of $X$ by the sets $X_s:=\{x\in X: f(x)\leq
s\}$. {\em Persistent homology} corresponding to such filtration is
the collection of morphisms 
\[H_k(m_{st}):H_k(X_s)\to H_k(X_t), k=0,1,\ldots
\]
induced by the 
natural inclusions 
$m_{st}:X_s\hookrightarrow X_t$ defined
for any $s<t$. (All homologies here are singular, over a field.)

The persistent homology theory (see, e.g. \cite{EH,
  Ou}) stipulates that for any $k\geq 0$, this collection of the
homology groups and homomorphisms between them 
decomposes into a direct sum (possibly, infinite) of elementary pieces, each of which is a trivial bundle of one-dimensional spaces over an interval, with trivialization being the isomorphisms. This means that a collection of homology classes in $H_k(X_b)$ emerges at some $s=b$ and ends at some $s=d$, $b<d$. (One augments this family of isomorphisms by 0-dimensional spaces outside of the $[b,d]$ range and corresponding zero morphisms $H_k(m_{st})=0$ whenever $s<b$ or $s>d$.) These chains of isomorphisms are referred to as {\em bars} $[b,d]$. 

A conventional and convenient way to visualize these bars is to place a point charge at
$(b,d)\in\Real^2$ in the half-space $b<d$ for the bar starting at $b$ and ending at $d$,
the weight being the rank of the corresponding homology group. This
results, for each $k\geq 0$, in what is known as the {\em
  $k$-dimensional persistent homology diagram $\ph_k$}.

We will identify persistence diagrams with their {\em counting measures}, i.e. the sum of
delta-functions weighted by the ranks of the homologies at $(b,d)$ for each bar. This gives a
positive measure (still denoted as $\ph_k$) supported by the
half-plane $\{b < d\}$.

For random functions, the resulting random measure is interpreted as a {\em point process}, a
convenient device to describe properties of the random persistence diagrams.

In general, the total mass of the persistence measures $\ph_k$ can be
infinite. (In fact, as we prove in \cite{BaWe}, it {\em is}
infinite for {\em generic} functions with loose enough moduli of continuity on triangulable
spaces with length metric.)

However, the contents of the displaced quadrants
$\quadr(b_*,d_*):=\{(b,d):b\leq b_*, d\geq d_*\}$ are finite for all
$b_*<d_*$ in most interesting situations. Filtrations (or functions generating them) for which this holds are called {\em tame} \cite{Ou}. Continuous functions on compact metric spaces are tame.

Notice that the points in $\quadr(b_*,d_*)$ correspond to bars {\em
  straddling} the interval $[b_*,d_*]$.

\subsection{Merge trees}

In this paper we deal exclusively with $0$-dimensional persistence homologies $\ph_0$.

It is well known that $0$-dimensional persistence diagram $\ph_0$ of the
filtration associated with a function $f$ can be defined without
invoking any homologies, but rather by simply tracking the connected components of
the subgraphs of the function. This procedure is codified by the
notion of {\em merge trees} associated with $f$ (our trees grow downwards).

Let $X$ be a path connected space \cite{burago2022course},  and $f$ a continuous real-valued function on $X$.

\begin{dfn}
For a pair of points $x,x'\in X$, we say that the point $x'$ {\em tops} $x$ 
if there is a path in $X$ from $x$ to $x'$, such that $f$ restricted to this path has a global maximum at its $x'$ end (clearly, this is a transitive relation). A pair of points topping each other are declared equivalent.

  The {\em Merge Tree} $\mtree_f$ of $f$ is the topological space which coincides as a set with the collection of these equivalence classes (we will denote the class of $x\in X$ as $[x]$).

It is clear that the function $f$ is constant over these equivalence classes, thus defining a function on $\mtree_f$ (which we still will be calling $f$). We will be referring to the (well-defined) value of $f$ at each point of the merge tree as its {\em height}.

The merge tree is equipped with the roughest topology such that the height is continuous. In this topology, the merge tree is a path metric space.\footnote{The metric trees, and their properties were used by probabilists with great success, see e.g., the survey \cite{legall}.}
\end{dfn}

Quite immediately, the merge tree of the merge tree $\mtree_f$ (viewed as a metric space with the height considered as the defining function) is again $\mtree_f$.

For compact spaces, the (unique) class of points where $f$ attains its global maximum on $X$ is referred to as the {\em root}. We note that one can extend the height of the root above the maximum of $f$ on $X$ (by attaching a long whisker sticking up);` sometimes it is useful to extend it to $+\infty$.

The class of a point that tops only points in its equivalence class is a {\em leaf} of the merge tree (leaves correspond to the local minima of $f$).

It is easy to see that classes of points topping a given point form a segment in $\mtree_f$, ending at the root: for any pair of levels $s<t$, and $x$ such that $f(x)=s$, there is a unique point $[x']\in \mtree_f$ topping $[x]$ with $f([x'])=t$.

This correspondence is, of course, the same as the natural mapping $H_0(X_s)\to H_0(X_t)$ used in the definition of the persistent homology: for each (path-)connected component of $X_s$, there is a unique connected component of $X_t$ containing it.

The term {\em merge} reflects the obvious from the definition observation that the components of the subgraph can only merge as the level grows, but cannot branch.

In general the merge trees can be quite wild even for Lipschitz functions (can have vertices of infinite degrees etc).


\subsubsection{From merge trees to $\ph_0$}\label{sss:mttoph}
Persistence diagrams can be derived from the merge tree by using the following recursive procedure (requiring compactness of sublevel sets) of {\em tree pruning}, relying on the {\em Elder Rule}\cite{EH}:
\begin{itemize}
  \item Find a global minimum (a leaf of the merge tree with the lowest height), and the unique path to the root of the merge tree from that leaf. The resulting pair of heights, from the lowest to the highest, is a recorded as a unit charge on the persistence diagram.
  \item We will refer to the path from the selected bottom leaf to the root as the {\em stem} of the (sub)tree. Removing the stem from the tree leaves a forest (possibly, empty) of merge trees. We will say that the stems of the each of the resulting subtrees {\em are attached to}, or {\em are descendants of} the stem just removed.
  \item Iterate the procedure on each of the resulting trees recursively. This results in a pile of stems, each corresponding to a bar in the persistence diagram.
\end{itemize}

From the construction of the stems, it is clear that a stem $[b',d']$ attached to its parent stem $[b,d]$ satisfies $b<b'<d'<d$.\footnote{We assume here, to avoid unnecessary in our context disambiguation efforts, that all critical values are different.}
In terms of the persistent diagrams, a descendent stem is located South-East off its parent.

The decomposition of a tree into a pile of stems erases a lot of information about the original function: thus, the parentage relationships between stems are lost. This is one source of ambiguity if one attempts to reconstruct the merge tree from its persistence diagram (pile of stems).

Another source of ambiguity (for univariate functions, which are the focus of this paper) is the fact that a stem can be attached to its parent on the right or on the left side. 

In combinatorial terms, the merge tree associated to a function on an interval is equipped the structure of a {\em plane tree}: this means that at each internal vertex one fixes the order in which branches are attached to it.

\subsubsection{Univariate Functions and Trees}
The merge tree (together with its planar embedding) generates in the standard fashion its {\em contour} or  {\em Dyck path} \footnote{One can define it, at least for the trees of the finite total length, as traversing the tree (imagined as a maze map) along the right wall in such a way that the distance to the root changes at the rate $\pm 1$; then the merge tree of the function equal to $-$ distance to the root will be equal to the original tree: see e.g., \cite{legall} for the standard definitions.}. 

In the situation when the underlying space is an oriented interval, the correspondence between height walks and plane trees is essentially bijective \footnote{In the standard convention, the contour walk starts at zero and remains positive; it is immediate how to reformulate it so that the heights of the vertices coincide with the corresponding critical values.}. 
  \begin{figure}[ht!]
    \begin{center}
      \captionsetup{font={small, it},width=.8\textwidth}
      \includegraphics[height=1.4in]{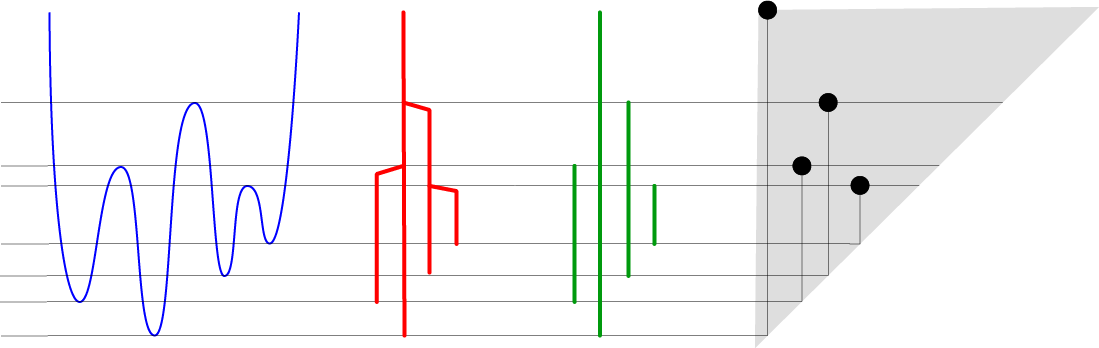}
      \caption{ Left to right: graph of a function, its merge tree, the pile of stems (bars) resulting from pruning, and the persistence diagram.}
      \label{fig:dyck}
    \end{center}
  \end{figure}


\begin{lem}
  The contour walk corresponding to the merge tree $\mtree_f$ of a continuous function on an interval is right-equivalent (i.e. equal up to a reparametrization of the argument) to $f$.
\end{lem}
The proof is immediate for the case of the finite trees: in this case the contour walk is a piece-wise linear function, with slopes $\pm 1$, interpolating between a finite sequence of critical values. It is immediate that such a piece-wise linear function is right-equivalent to the original function (reparameterizing the timeline so that on each interval of monotonicity of the original function it becomes linear after the coordinate change). For general continuous functions, the proof is more involved and won't be presented here, as we do not use this Lemma in the sequel.

Figure \ref{fig:dyck} illustrates the relationship between univariate functions, their merge trees and their $\ph_0$.

\subsubsection{From $\ph_0$ to merge trees}
Constructions above essentially answer the natural question about the space of univariate functions generating a particular collection of bars as its $\ph_0$ diagram. 
The extra data necessary to rebuild the planar merge tree (and therefore the function generating that merge tree, up to reparametrization) from a (locally finite) collection of bars, are the data a) on their parental attachments, and b) on their planar order.

This result was first published by J. Curry \cite{C}: Let $f$ be a univariate function with the longest bar $-\infty\leq b_\infty<d_\infty\leq \infty$. Assume that the number of bars straddling each interval (i.e. the $\ph_0$-content of $\quadr(b,d)$) is finite (i.e. the function is tame).

\begin{prp}[ \cite{C}]\label{prp:curry}
  The merge tree is uniquely reconstructed by the (arbitrary) attachment relation on the stems (how each bar with exception of the longest one is attached to one of the straddling bars), and by the planar structure (which we will refer to as chiralities, see section \ref{subsec:Chir}) of the stems.
\end{prp}

We remark that one can generalize Proposition \ref{prp:curry} to higher dimensions as follows:

\begin{thm}
Consider the space of Morse functions on $\Real^d$ with compact sublevel sets, whose critical points are of indices $0$ and $1$ only, and increasing to infinity at infinity. The subspace of functions corresponding to a given merge tree (i.e. with a fixed persistence diagram and the stem attachments data) is homotopy equivalent to the product of $(d-1)$-dimensional spheres, one sphere for each bar not ending at $+\infty$.
\end{thm}

The proof will appear elsewhere.
\subsection{Chirality}\label{subsec:Chir}
A bar (i.e. a point $(b,d)$ in the persistence diagram) corresponds to two events: the birth of a connected component of the sublevel set $\{f<h\}\}$ at height $h=b$, and its death, i.e. merger of two connected components, at height $d$. Clearly, this implies that there exists a (local) minimum of $f$ with the critical value $b$ and a (local) maximum with the critical value $d$.

\begin{dfn}\label{dfn:coupled}
We will be referring to the critical {\em points} corresponding to a bar $[b,d]$ in $\ph_0$ as {\em coupled}.
\end{dfn}

Given a coupled pair of local maximum and minimum, we will refer to the right or left orientation of $0$-dimensional bar in a persistence diagram of a univariate function, when it is attached to its parent stem, as {\em chirality}. There are two possible chiralities which we will denote as \maa, and \fra.\footnote{\ma is the 11th letter of Ukrainian alphabet, phonetically close to {\tt /y/}.} In terms of the bijection between the Dyck paths and plane trees, \ma corresponds to the branch attached to its parent stem from the left; \fra, from the right.

Formally,
\begin{dfn}
  A bar $(b,d)=(f(s), f(t))$ (here $s,t$ are the critical points corresponding to the critical values
  $b,d$) is an \maa, if
  the local maximum follows local minimum, i.e. if $s<t$, and an \fra, otherwise, - i.e. if $s>t$.
\end{dfn}
  \begin{figure}[ht!]
    \captionsetup{font={small, it},width=.8\textwidth}
    \begin{center}
      \includegraphics[height=1.in]{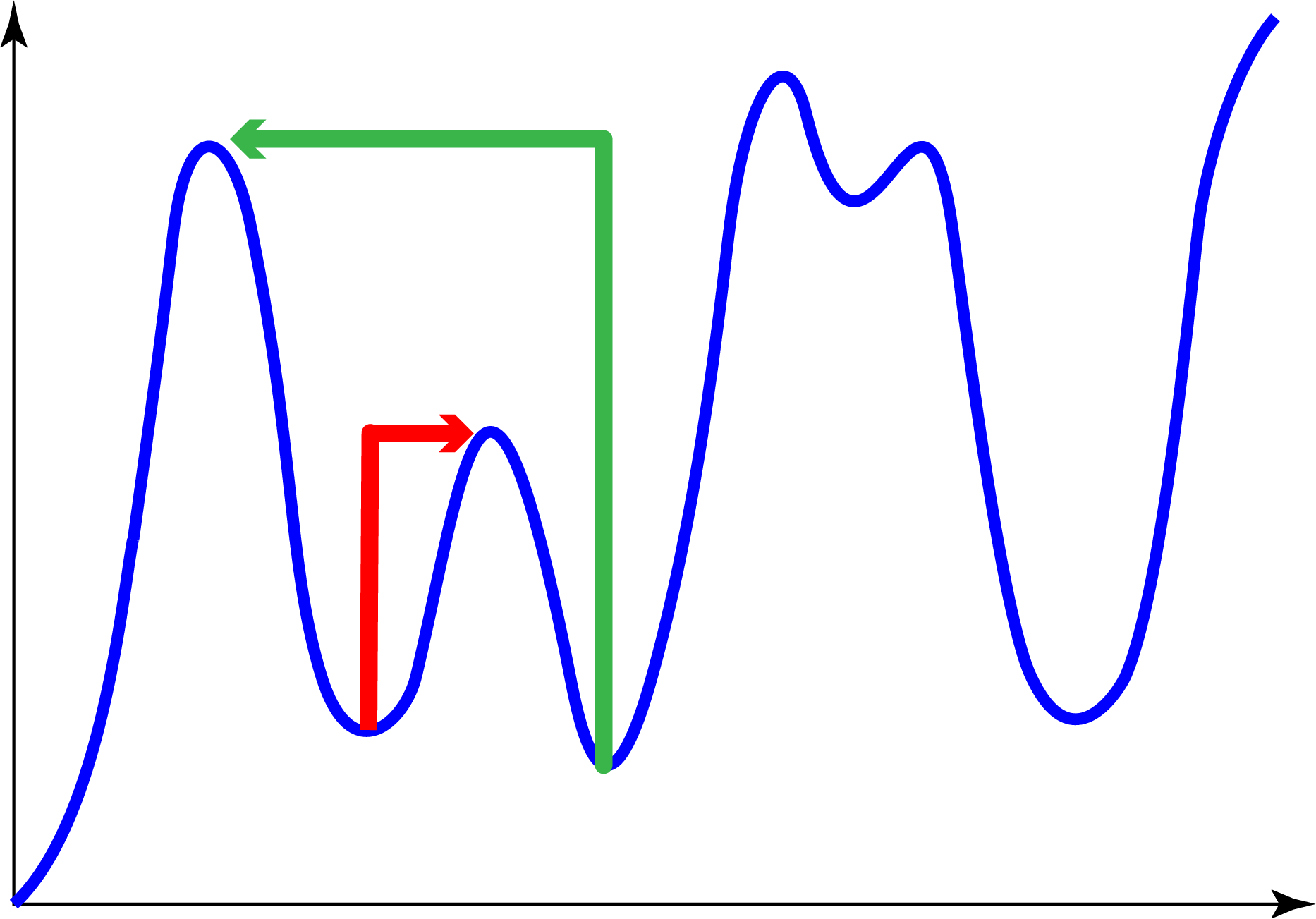}
      \caption{Coupled critical points: \ma in red and \fr in green.}
      \label{fig:chirality}
    \end{center}
  \end{figure}
\begin{rem} Chiralities might be a useful tool to capture potential asymmetry of a time series with respect to the time reversal, a prominent topic in econometric literature.  In \cite{lacasa_time_2012}, an approach somewhat resembling our chiralities was proposed.
\end{rem}

Intuitively, if the function is rising, one would expect more of \fra's than \maa's. To see this, recall that for a univariate function $f$, a local minimum $x$ is coupled by a bar to a local maximum $y$ (with critical value $h=f(y)$), if 
a) the connected component of the sublevel set $f<h$ containing $x$ has that point as the global minimum, while b) the global minimum of the other component sits deeper than $f(x)$.
What happens if a small linear tilt is added, so that $f(t)\mapsto f(t)+s(t-y), s>0$ (we can assume that this linear function vanishes at the local maximum as adding constants does not affect persistence couplings)?

If $y<x$ (i.e., we have a \fr chirality), the minimum of $f$ over the right component will rise, while that over the left component will go down, and the coupling will remain unchanged. If, on teh other hand, $x>y$ (i.e., we have a \ma chirality), the minimum of $f$ over the left component can become deeper than the minimum of $f$ over the right one, and the coupling can flip.

Hence, if started with a random function invariant with respect to time reversal, one would expect that a tilt with positive slope would generate more \fra's than \maa's.

In Section \ref{sec:chir_brown} we will make this intuition precise for Brownian trajectories with positive drift, by computing the average excess of the number of \fra's  over \maa's.

\subsection{Bars and Windings}\label{sec:bars}

In this section we will set up necessary apparatus to analyze bar decomposition of a continuous univariate function. Our key observation is that bars in a persistence diagram for a function on the real line are in direct correspondence with the {\em windings of the function around an interval}.\footnote{The term is clarified on Figure \ref{fig:auto_a}.}

In what follows we will be working primarily with the functions on an interval which attain their global minimum on the left end, and their global maximum on the right. One can always relate this particular version with other settings (say, where the function tends to its global supremum on both ends), by slightly expanding the domain of the univariate function $f:[t_0,t_1]\to\Real$ to $[t_0-c, t_1+c], c>0$ and defining $f$ on $[t_0-c,t_0]$ to be a linear function interpolating between $\min_{[t_0,t_1]}f-1$ and $f(t_0)$, and on $[t_1,t_1+c]$ to be a linear function interpolating between $f(t_1)$ and $\max_{[t_0,t_1]}f+1$. Such a transformation alters the persistence diagram by adding at most by one (long) bar.

\subsubsection{Winding of a function around an interval}\label{ss:windings}
Recall that we say that a point $(b,d), b<d$ in the persistence diagram (or the corresponding bar $[b,d]$) {\em straddles the interval $[b',d']$} if
$b\leq b'<d'\leq d$.

Given an interval $[b,d]$, and a continuous function $f$ on $[t_0,t_1]=I\subset\Real$, starting below $b$ and ending above $d$ on $I$ (i.e. $f(t_0)<b<d<f(t_1)$), one can define sequences of alternating $b$- and $d$-times $t_0^b<t_1^d<\ldots<t_k^d<t_k^b<\ldots$ by setting, iteratively,
\begin{equation}\label{eq:windings}
  t_0^b=t_0; \mathrm{\ for\ } k\geq 1, t^d_{k}=\min\{t:f(t)=d, t\in I, t>t^b_{k-1}\}, t^b_{k}=\min\{t:f(t)=b, t\in I, t>t^d_k\}.
  \end{equation}
  (The iterations stop if the set over which the last $\min$ is defined is empty.)

This sequence is obviously finite, as the interval $I$ is compact and $f$ is continuous.

  \begin{dfn}\label{def:windings}
    A {\em winding of $f$ around the interval $[b,d]$} is a pair $\{t_k^b,t_k^d, k\geq 1\}$ in the sequence thus generated.
  \end{dfn}
  (See Fig. \ref{fig:auto_a}, right display, for an example of a function with two windings around an interval.)
  We remark, that the continuity of $f$ implies that the set of windings in any compact subset of $I$ is finite.

  One has the following
  \begin{lem}
    The number of the bars in $\ph_0(f)$ straddling the interval $[b, d]$, or, equivalently, the $\ph_0$ content of the quadrant $\quadr(b,d)$, is equal to the number of windings of $f$ around $[b,d]$, plus one.
  \end{lem}
  \begin{proof}
    On one hand, to each winding $(t^b_k, t^d_k)$ one can associate the unique connected component of $I_d:=\{f<d\}$ containing $t_k^b>t_0$: each such component generates a bar. On the other hand, to any component of $[t_-, t_+]\subset I_d$ over which $\min_{[t_-, t_+]} f<b$ (besides the one having $t_0$ as its left end), one can associate the winding with the left end $t^d_k=t_-$. This is manifestly a bijection between the windings and all but the leftmost components of $I_d$ over which $\min f\leq b$.
  \end{proof}
  This definition of windings is, of course, identical to the well-known construction of "upcrossings" used in the standard proof of Doob's martingale convergence theorem \cite{doob}. 

    \section{Persistence Diagrams for Brownian Motions with Drift}\label{sec:brownian} 
    Our focus in this note on the properties of persistent diagrams for Brownian motions, the simplest Markov processes with continuous trajectories.

    Specifically, we consider Brownian motion with a drift, that is
    \be
    f(t)=\brm_o(t)+mt,
    \ee
    where $\brm$ is the standard Brownian motion starting at $0$, on the ray $(0,\infty)$. While our definitions in the preceding sections handled only finite intervals, the Brownian motions with positive drift will never return to any upper-bounded set (such as $(-\infty,d]$) after a finite random time almost surely, allowing us to use these constructions without modifications.
    
    We will be dealing with the intervals $0<b<d$. For general intervals $[b,d]$ with $b<0$ we can condition on the process $f$ to ever reach the level $b$ (this happens with probability $\exp(-2mb)$); restarting the process at the stopping time of the first hitting $b<0$ allows one to immediately generalize the results below to the case $b<0$.
    \subsection{Persistence Diagrams as Point Processes}
    Almost surely, there is unique bar $(b_\infty,d_\infty)=(\inf_{t>0} f(t),\infty)$.

    Other bars are finite, forming a random persistence diagram $\ph_0$. Our interpretation of this persistence diagram is to view it as a {\em point process}, a random measure consisting of sum of $\delta_{(b,d)}$ over all bars $(b,d)$.

 \begin{rem} Almost sure tameness of the this point process, that is the fact that the content of any NW-quadrant with its apex above the diagonal is almost surely finite, will be proven below.     
However, this can be seen from the mentioned above general result in \cite{BaWe}, that the persistence diagram of any H\"older function on a compact metric space is tame, and the well-known fact that the trajectories of Brownian motions (with any smooth drift) are almost surely $c$-H\"older for any $c<1/2$.
   \end{rem}

    \subsubsection{Process decomposition}
    Following the definitions of Section \ref{ss:windings}, we consider the alternating sequences of $b$- and $d$-times. 
    
    These times are, clearly, stopping times with respect to the filtration
    $\{\sif_\tau\}_{\tau\geq 0}$ 
    generated by the Brownian motion. Hence, the intervals $t^b_k-t^b_{k-1}, k=1,2,\ldots$ between the winding events are independent (thanks to the standard fact that for any stopping time $\tau$, the process $\brmd{m}_{\tau+t}-\brmd{m}(\tau)$ is independent of $\sif_\tau$).
    As an immediate  consequence, we obtain

    \begin{prp}\label{cor:geom}
      Consider the standard Brownian motion with constant drift $m>0$.
      Then the total number of bars in $\ph_0$-persistence diagram straddling an interval $[b,d], b>0$ (or, equivalently, the $\ph_0$-content of the quadrant $Q(b,d)=\{(x,y):x<b<d<y\}$)
      is geometrically distributed with parameter $p=p_m(\Delta)=1-\exp(-2m\Delta)$, where $\Delta=d-b$.
    \end{prp}

    \begin{proof} Indeed, the total number of windings around $[x,y]$ is the total number of  transitions from $d$ to $b$. The probability that the trajectory of the standard Brownian motion with drift $m$ starting at $y$ reaches $x$ before escaping to infinity is $\exp(-2m(y-x))$. Strong Markov property implies that the escapes to infinity after each hitting of $y$ are independent, proving the conclusion.
    \end{proof}

    \subsubsection{Invariance properties}\label{sec:invar}
    Let us state explicitly two invariances.

    First, the standard rescaling properties and the invariance  of the persistence diagrams under the reparameterization of the argument imply that the same result - and all results on the distributional properties of the persistence diagrams and chiralities for the standard Brownian motion and drift $m$ remain valid for the Brownian motion with quadratic variation $\sigma^2$ and drift $m\sigma^2$.

    Further, for $b\geq 0$ the distribution of the $\ph_0$ content of $\quadr(b,d)$ depends only on the length $\Delta$ of the bar $[b,d]$. This is again an immediate property of the Markov property of the process. Similarly, the joint distributions of the windings around several intervals remains invariant under the simultaneous shifts of these intervals as long they all remain within the positive half-line.

    \subsubsection{Other Brownian motions}\label{sss:otherbm}
    We should notice that reflection principle allows one to derive quite easily the expected number of bars straddling a particular interval for other Brownian motions. For example, for the standard Brownian bridge on the unit interval, and an interval of $[b,d], 0<b<d$, it is equal to
    \[
    \ph_0(\quadr(b,d))=\sum_{n\geq 1} \exp(-2(n\Delta+b)^2),
    \]
(with $\Delta=d-b$). For the intervals containing the origin, the expression is a bit more cumbersome, but still easy to derive.

    \subsubsection{Intensity measure for $\ph_0$ for Brownian motion with drift}
    Corollary \ref{cor:geom} immediately leads to an explicit expression for the intensity
    measure $\ex\ph_0$ for the point process of zero-dimensional persistence from standard Brownian motion with drift:

    \begin{prp}\label{prp:intensity}
      The intensity measure of $\ph_0$ on $0<b<d$ is
      given by the density
      \be\label{eq:intensity}
      \mph=\frac{4m^2e^{-2m\Delta}(1+e^{-2m\Delta})}{(1-e^{-2m\Delta})^3},
      \ee
      where $\Delta=d-b$.

      In particular, near the diagonal the density explodes as $1/m(d-b)^3$.
    \end{prp}

    \begin{proof}
      By \ref{cor:geom}, the expected $\ph_0$ content of the quadrants $Q_{b,d}$ is
      \[
      q_{b,d}:=\ex \ph_0(Q_{b,d})=p_m(\Delta)/q_m(\Delta)=\frac{1}{e^{2m\Delta}-1}.
      \]
      The density of the intensity measure is, clearly, given by
      \[
      \mph(b,d)=-\frac{\partial^2 q_{b,d}}{\partial b\partial d},
      \]
      which gives the stated result.
    \end{proof}

The fact that the density grows as $(m\Delta)^{-3}$ as $\Delta\to 0$ is consistent with the scaling laws of \cite{BaWe}, and holds for other Brownian motions (for Brownian bridges this follows from the formula in section \ref{sss:otherbm}), see also \cite{perez}.

    \section{Chiralities in Brownian trajectories}\label{sec:chir_brown}
    Fix an interval $\vertical=[b,d]$. Corollary \ref{cor:geom} describes the distribution of the total number of the bars $\ph_0(Q(b,d))$ of the Brownian motion with drift $\brmd{m}$ straddling $\vertical$. In this section, we address the question: {\em how many of them will be \maa,  and how many \fra}?

    \subsection{Deconstructing into the bars}
    Define, as in \eqref{eq:windings} the stopping times $t_l^b, t_l^d, l=1,\ldots$. A trajectory straddles the interval exactly $k$ times if $(k+1)$ is the first index for which $t_{k+1}^b=\infty$. In this case there are $k$ pairs of times, $(t_l^d< t_l^b)$ forming an increasing sequence collectively: $t_0^b<t^d_1<t^b_1<t^d_2<t^b_2<\ldots<t^d_k<t^b_k<t^d_{k+1}$ (it is convenient to add the boundary points to the sequence).

    We will use these times to chop the trajectory $\brmd{m}$ into $2k$ pieces, $f^\downarrow_l, f^\uparrow_l$ by setting
    $$
    f^\downarrow_l(t):=\brmd{m}(t^d_l+t), 0\leq t\leq t^b_{l}-t^d_l, 1\leq l\leq k
    $$
    being the fragments of the trajectory traveling from $d$ to $b$, and
    $$
    f^\uparrow_l(t):=\brmd{m}(t^b_l+t), 0\leq t\leq t^d_{l+1}-t^b_l, 1\leq l\leq k
    $$
    their upward counterparts.

    \begin{prp}\label{prp:chop}
      Conditioned on the event that the trajectory straddles $[b,d]$ $m\geq k$ times, the random processes $f^\uparrow_l, f^\downarrow_l, l=1,\ldots,k$ are independent. Moreover the processes $f^\uparrow_l, l=1,\ldots,k$ are identically distributed, as the processes $f^\downarrow_l, l=1,\ldots,k$ are.
    \end{prp}
    \begin{proof} This follows, again, from the strong Markov property of the Brownian motion with drift.\end{proof}

      We remark that the distributions of the downward segments $f^\downarrow$ are those of the standard (driftless) Brownian motion starting at $d$ and stopped once it reaches $b$, while the distributions of $f^\uparrow$ are those of the Brownian motion with drift $m$ started at $b$ and stopped when it reaches $d$.

      The {\em set} of the maximal values of the processes $f^\downarrow_l, l=1,\ldots,k$ coincides with the set of the right ends of the bars straddling the interval $[b,d]$, while the minimal values of the processes $f^\uparrow_l, l=1,\ldots,k$ form the set of their left ends. Of course, $d_l$ is not necessarily coupled to $b_l$ (although it could be).

      \subsection{Symmetries of the process}
Almost surely, these local maximal and minimal values are distinct. When this is the case, we form two {\em rank} processes, $\sigma_+(l)$ and $\sigma_-(l)$, where $1\leq l\leq k$, $\sigma_+$ is the rank of $d_l$ among all $\{d_1,\ldots, d_k\}$ and $\sigma_-(l)$ is the rank of $b_l$ in $\{b_1,\ldots, b_k\}$. The sequences $\sigma_+$ and $\sigma_-$ are permutations of $\{1,\ldots,k\}$.

      The Proposition \ref{prp:chop} implies

      \begin{prp}
        Conditioned on the $\ph_0(Q_{a,b})= k$ the permutations $\sigma_+$ and $\sigma_-$ are independent,
        uniformly distributed on the symmetric group $\mathfrak{S}_k$.
      \end{prp}
      \begin{proof}
        Both statements follow immediately from the mutual independences of $f^\uparrow_l, f^\downarrow_l,\  l=1,\ldots,k$.
      \end{proof}

\begin{figure}[htbp]
        \captionsetup{font={small, it},width=.8\textwidth}
        \begin{center}
          \includegraphics[height=1.8in]{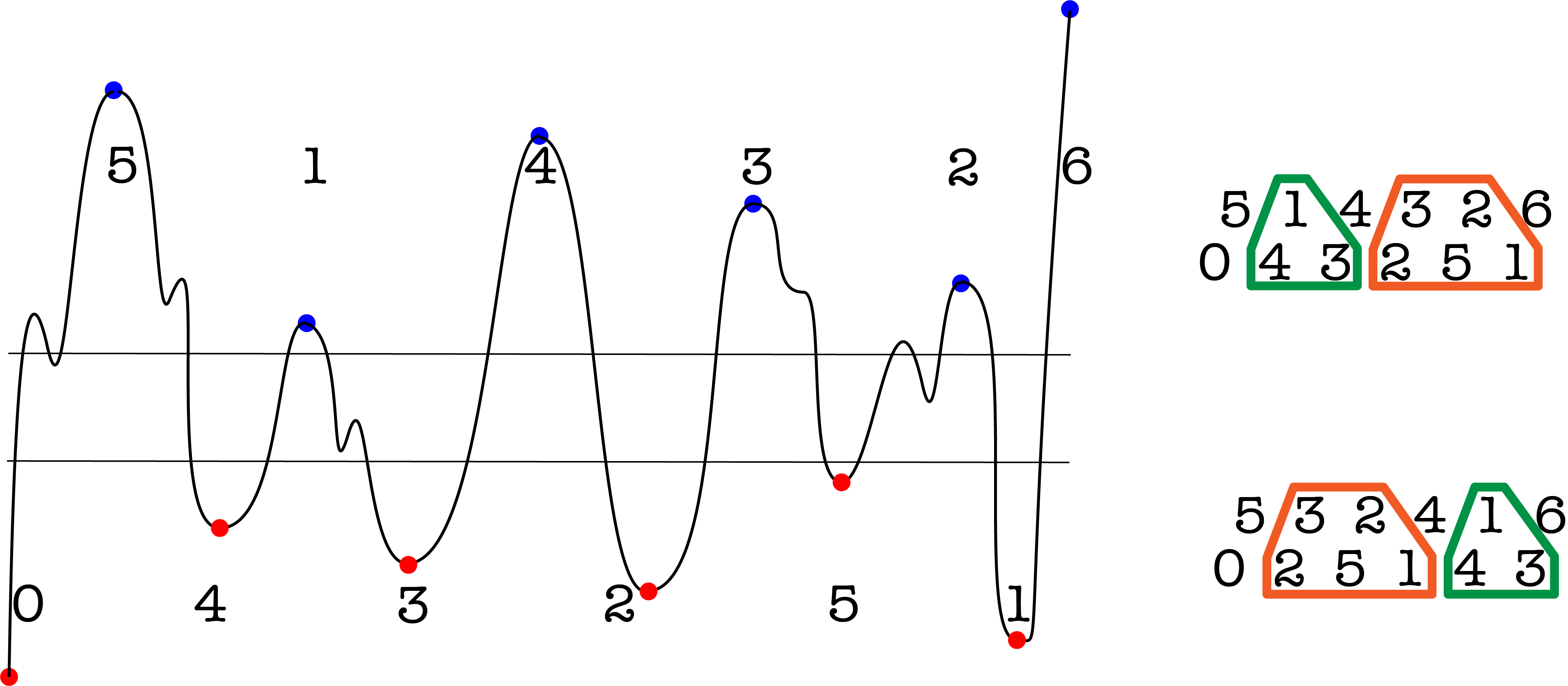}
          \caption{{\em Left display:} There are $5$ bars straddling $I$ in this picture. The permutations of maxima and minima are, correspondingly, $\sigma_+=(514326)$ and $\sigma_-=(0 43251)$. The walls for $r=\sigma_+(3)=4$  are $s=5$ and $t=6$, at the positions $1$ and $6$ respectively, the left well is $3$, the right well is $1$ and $r=4$ is coupled with the $\sigma_-(2)=3$, thus a \maa.
{\em Right display:} swapping the segments.}
          \label{fig:walls}
        \end{center}
      \end{figure}

      Which of the local minima and maxima $b, d$ are coupled, depends only on the permutations $\sigma_+, \sigma_-$. The left display of Fig. \ref{fig:walls} shows a configuration of maxima and minima of the processes $f^\uparrow_\cdot, f^\downarrow_\cdot$. 
      
We will augment the sequence $\sigma_+$ with $(k+1)$ on its right, and the sequence $\sigma_-$ with $0$ on its left.

For a maximum of rank $r, 1\leq r\leq k$ at position $l$ we call the rightmost index $m<l$ for which $\sigma_+(m)>\sigma_+(k)$ its {\em left wall}. Note that the left wall does not exist if (and only if) $\sigma_+(l)$ is a {\em record} in the sequences $\sigma_+$ (that is it beats all the earlier elements). Thus $\sigma_+(1)$ is always a record.

Similarly we define the {\em right wall} of $\sigma_+(l)=r$, as the leftmost index $n>l$ for which $\sigma_+(n)>r$. Note that thanks to our convention that $\sigma_+(k+1)=k+1$, the right wall always exists.

We say that the maximum of rank $r$ is walled by ranks $s,t>r$ if $s$ and $t$ are ranks of the left and right walls of the rank $r$ maximum (in particular, implying that $r$ is not a record in the permutation $\sigma_+$). In this situation, we denote the positions of the elements $s,r,t$ as $m,l,n$.

Summarizing, for a two-side walled (i.e., non-record) element $r$,
\be
\sigma_+(m)=s>\sigma_+(l)=r<\sigma_+(n)=t, m<l<n \mbox{ and } 
\sigma_+(u)<r \mbox{ for }m<u<l \mbox{ or } l<u<n. 
\ee
If element $r$ at position $l$ is walled by $s,t$ at positions $m,n$ we will form its {\em left} and {\em right} wells, defined as 
\be
\mathtt{LV}(r):=\min_{m\leq u<l}\sigma_-(u);\mathtt{RV}(r):=\min_{l\leq v<n}\sigma_-(v).  
\ee

\begin{lem}
If element $r$ at position $l$ is walled by $s,t$ at positions $m,n$ then it is coupled with one of its wells: with the left well if $\mathtt{LV}(r)>\mathtt{RV}(r)$ and the right well otherwise.
\end{lem}

\begin{proof}
The components merging at $r$ are born, obviously, at its wells. The choice of the well to couple with follows from the Elder Rule.
\end{proof}

Thus, if the left well is deeper, the bar corresponding to $r$ couples to the right (i.e. is an \fra); if the right one is deeper, it is an \maa.

Conforming with the intuition, if the element is a record in $\sigma_+$, its left well at $0$ is the deeper one, and the element is coupled to an index on its right, that is, it is a \fra.

\subsection{More \fra's than \maa's}
      Finally we can deduce the excess of  the expected numbers of \fra's over \maa's in a Brownian motion with drift.

      \begin{thm}\label{thm:foverm}
        Conditioned on the number $k$ of bars straddling an interval $I$, in a trajectory of $\brmd{m}, m>0$, the expected excess of the numbers of \fra's over \maa's is equal to the expected number of records in a random $k$-permutation, that is $k$-th harmonic sum,
        \be\label{eq:harmonic}
        \ex(|\fra|-|\maa|)=H_k=1+1/2+\ldots+1/k.
        \ee
      \end{thm}
      \begin{proof}
      The couplings of the elements of $\sigma_+$ to the elements of $\sigma_-$ depend only on those permutations, so we can focus on them.
      
      Fix a rank $r$. Partition the set of pairs $\Sigma:=(\sigma_+, \sigma_-)$ into $\Sigma_*$, corresponding to permutations $\sigma_+$ for which $r$ is a record, and 
      \[
      \amalg_{s,t>r, s\neq t}\Sigma_{s,t}
      \]
      of permutations for which $s,t$ are respectively the left and the right walls for $r$.
      
      One can define an involution on each of the collections $\Sigma_{s,t}$. If the position of $s, r, t$ are, respectively, $m<l<n$, then swap the segments 
\[
(\sigma_+(m+1),\ldots\sigma_+(l-1)) \mbox{ and } (\sigma_+(l+1),\ldots\sigma_+(n-1))
\]      
(keeping $\sigma_+(l)=r$ between them), and swap the segments
\[
(\sigma_-(m),\ldots\sigma_-(l-1)) \mbox{ and } (\sigma_-(l),\ldots\sigma_+(n-1))
\]
(right display on Fig.\ref{fig:walls} shows this swap for a configuration in $\Sigma_{5,6}$ for $r=4$).

Under such a swap, the chirality of the coupling changes sign (as the wells swap). This shows that the contributions from the partitions $\Sigma_{s,t}$ to the difference between the numbers of \fra's and \maa's coupled to $r$ over these subsets vanishes.

Thus only the records contribute. The probability that $r$ is a record is $1/(k-r+1)$ (or, see, e.g. \cite{pitman_combinatorial_2006}), and the result follows. 
      \end{proof}

      \begin{cor}
        For  Brownian motion with drift $m$ , the expected excess $\ex(|\fra|-|\maa|)$ in the bars straddling an interval of length $\Delta$ is
        $$
        -\log(1-\exp(-2m\Delta)).
        $$
      \end{cor}
      \begin{proof}
        The sum \eqref{eq:harmonic} gives the excess conditioned on the number of straddling bars. Summing these harmonic sums with weights $pq^k$, where $q=\exp(-2m\Delta)$, the probability given by Corollary \ref{cor:geom}, gives the desired result, after switching order of summation in the (obviously) converging series.
      \end{proof}

      In particular, the expected excess of \fr\  over \ma\  grows logarithmically (as $|\log(2m\Delta)|$) for small $m\Delta$; consequently the fraction of the excess among all bars straddling a short interval, disappears, as the length of the interval decreases to $0$.

      \section{Automata and Correlation Functions for $\ph_0$.}\label{sec:covar}
      In this section we will deploy somewhat more involved techniques to understand the correlations between the windings around different intervals (which are essential to understanding the correlation functions of the $\ph_0$ point process). This relies on the constructions of finite state automata translating continuous functions into symbolic trajectories. For continuous trajectories given by strong Markov processes, thus generated symbolic trajectories become Markov processes, and the questions about the joint distributions of windings around different intervals become standard questions about numbers of specific transitions in realizations of these Markov chains. We will use these ideas to derive some second-order characteristics of the $\ph_0$ point processes.
      
Recall that any point process is characterized by its {\em factorial moments}, which in our tame setting can be characterized by their densities, i.e. {\em correlation functions} \cite{ruelle}.

The goal of this section is to outline a formalism to deduce those correlation function, and to compute a second order correlation function $\mph_2$ in some situations.
      
 \subsection{Windings and Symbolic Trajectories}\label{subs:automata}
  To study windings of a function around an interval, or more generally, the interactions of the function with one or more intervals, it is convenient to encode them with a scheme using finite automata. We will use an intermediate construction that will turn a continuous function of time into a symbolic trajectory.

  \subsubsection{Coverings and Lifts}\label{sss:lifts}
  Consider an open interval $U=(h^-,h^+)$ in the real line, and its finite covering by open intervals 
  \[
  U=\cup_{\inl\in\intervals} U_\inl,  U_\inl=(h^-_\inl,h^+_\inl), \inl\in\intervals.
  \]
  We will consider the disjoint union of these intervals,
  \[
  \tU:=\amalg_{\inl\in\intervals} U_\inl
  \]
  together with the natural projection, $\pi:\tU\to U$. While this is certainly not a covering, one can define {\em lifts} of continuous trajectories in $U$ to $\tU$, with a few bits of extra data.

  We will denote by $\boundaries$ the collection of {\em finite} ends of these intervals: 
  \[
  \boundaries:=\{h_\inl^\pm: h^-<h^\pm_\inl<h^+\}.
  \]
  Fix a map $\tau$ from this set $\boundaries$ of boundaries to the intervals $\intervals$, subject to the condition that each finite end $h^\pm_\inl$ is mapped into a (necessarily, different) interval which contains it:
  \[
\mbox{ if } \tau(h^\pm_\inl)=\inl', \mbox{ then }  h^\pm_\inl\in U_{\inl'} \mbox{ and, of course, } \inl'\neq \inl.
\]

Given these conditions, a continuous function $f:[t_s, t_f]\to U$, can be lifted unambiguously to a function $\tf:[t_s, t_f]\to \tU$ such that $\pi\circ\tf=f$, provided one assigns the initial interval $U_{\inl_s}$ to which $f(t_s)$ is lifted (so that $f(t_s)\in U_{\inl_s}$). 

Such a lift is given by a finite partition of the time interval $[t_s, t_f)=\amalg_k [t_k,t_{k+1})$, and a sequence of intervals $U_{\inl_k}$ and continuous functions $\tf_k:[t_k,t_{k+1})\to U_{\inl_k}$ defined inductively as follows:
      \begin{itemize}
      \item One initializes with $t_0=t_s, \inl_0=\inl_s$;
      \item For $t_k\leq t<t_{k+1}$, one has $\tf(t)\in U_{\inl_k}$ (satisfying $\pi\circ\tf=f$);
      \item The switch times $t_{k+1}$ are defined as $\inf\{t>t_k:f(t)=h^+_{\inl_k} \mbox{ or } h^-_{\inl_k} \}$, and 
      \item The next interval is given by the transition map $\tau$: $\inl_{k+1}=\tau(h_{\inl_k}^\pm)$ ($\pm$ depending on which end of $U_{\inl_k}$ the trajectory hits at $t_{k+1}$).
      \end{itemize}

      In other words, the $k$-sojourn lasts while the trajectory stays within the corresponding interval; as soon as it hits a boundary, it instantaneously jumps to whichever interval the switch map $\tau$ sends it, all the while tracking the function $f$.
      \begin{rem}
        This construction of lifting a continuous trajectory can be easily generalized to, say, coverings of manifolds by the open patches with well-behaving boundaries: one needs to partition those boundaries into disjoint subsets such that each of them is contained in a patch; the switch map would then send the trajectory hitting the boundary at a component to one of the patches containing that component. This construction is an alternative way to define the diffusions on cellular complexes, a notoriously messy problem (compare \cite{nye_random_2020}).
      \end{rem}

\subsubsection{Finite State Automata}\label{sss:automata}    
Given the data described in Section \ref{sss:lifts} we can form finite state automata whose states will correspond to the finite endpoints of intervals $U_\inl$. 
A bit counterintuitively, the state corresponding to $\inl^\pm$ will be located in the interval $\tau(\inl^\pm)$, i.e., at the {\em exit} point of the switch. Once located in that interval, the process can exits it only through one of the ends.
Therefore, in terms of the (directed) graph description of the automaton, each state will have at most two potential out-arrows.

We will fix an interval $(h^\alpha,h^\omega)$ containing all the intervals $U_\inl, \inl\in\intervals$, referring to it as the {\em span}.

We add the absorbing state $\omega$ and the initial state, $\alpha$ corresponding to the endpoints of the span. The initial state $\alpha$ has no in-, and the absorbing state $\omega$, no out-arrows.

\subsubsection{Symbolic trajectories}\label{sss:symbolic}
  Fix a finite automaton described in \ref{sss:automata} with the span $(h^\alpha,h^\omega)$.

  Consider a continuous function $f$ {\em crossing the span of the automaton}, and such that on its interval of the definition $I=(t_0,t_1)$ contains the lifespan: such a function starts below the span of the automaton, and ends above it (i.e. $f(t_0+)<h^\alpha; f(t_1-)>h^\omega$). 
  
Each function crossing the span defines a finite {\em symbolic trajectory}: a sequence of states, starting with $\alpha$, ending with $\omega$ and transitions between finite states in between. These transition track the exits of the lift of $f$ to the intervals $U_\inl$, as described in Section \ref{sss:lifts}.

  \begin{figure}[ht!]
    \captionsetup{font={small, it},width=.8\textwidth}
    \begin{center}
      \includegraphics[height=2in]{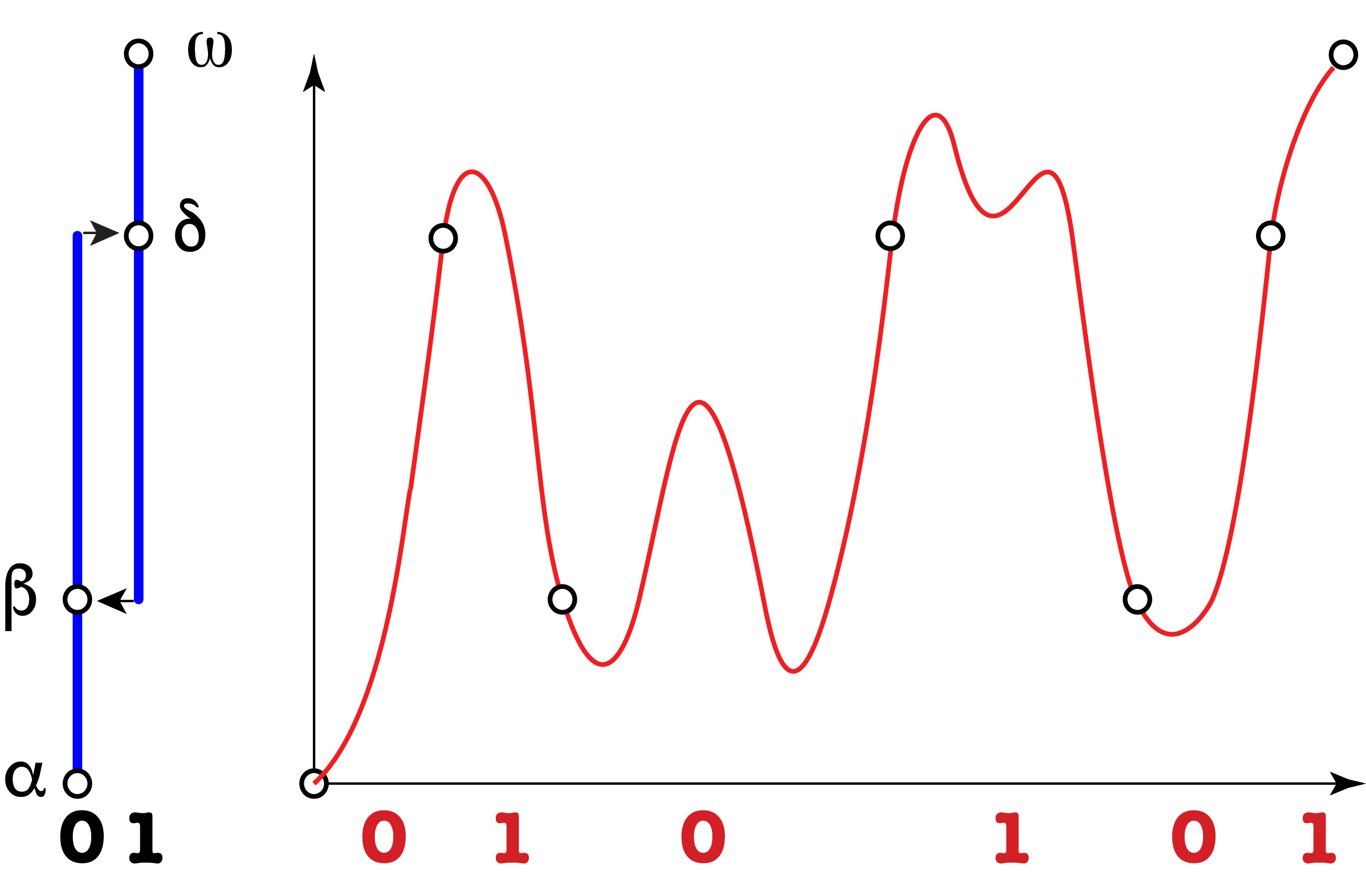}
      \caption{The finite automaton justifies the term {\em windings}: the number times a function straddles an interval is equal to the number of times its lift cycles around. The symbolic trajectory corresponding to the function shown above is $\alpha\to\delta\to\beta\to\delta\to\beta\to\delta\to\omega$. There are two transitions $\delta\to\beta$, and two windings around the interval $[b,d]$.}
      \label{fig:auto_a}
    \end{center}
  \end{figure}

  The resulting sequence of states and arrows is called {\em symbolic trajectory} corresponding to the function $f$.

  The following is immediate:
  \begin{lem}
    The symbolic trajectory of a continuous function crossing the span of a finite automaton is finite.
    \end{lem}\qed

    \subsubsection{Automata and windings around an interval}
    The simplest useful finite automaton of this kind described above (shown on the left display of Fig. \ref{fig:auto_a}) has the starting state $\alpha$ at level $h^\alpha$; two finite states $\beta$ and $\delta$ with levels $b:=h^\beta$ and $d:=h^\delta$, and the absorbing state $\omega$ with level $h^\omega$. The transitions are $\alpha\to\delta$, $\delta\to \beta, \delta\to\omega$, and $\beta\to\delta$.

    Again, the following is immediate:
    \begin{prp}\label{prp:number}
      For a continuous function $f$ crossing the span of this automaton, the number of $\delta\to\beta$ transitions in the corresponding symbolic trajectory is one less than to the number of bars straddling the interval $(b,d)$ in the persistence diagram $\ph_0$ of $f$.
    \end{prp}

    (Remark: The one extra bar is the "infinite" one, spanning the $\inf f$ and $\sup f$.)

\subsection{Lifts of Brownian Paths}
The formalism of lifts of trajectories to the automata described in Section \ref{sss:automata} is well suited to address the windings of Brownian motions around intervals.
Namely, in our model situation, --- where the lifted function is a continuous trajectory of a strong Markov process, --- the symbolic trajectory becomes a Markov process itself. 
      
\subsubsection{Transition Probabilities for Brownian Motion with Drift}

It is easy to compute the transition probabilities in the situation of Brownian motion with drift. 

Let's introduce the character defined by
      \be
      \Em(s)=\exp(-2ms).
      \ee
Evaluated on the Brownian motion with drift $m$ it becomes a martingale: hence its ubiquitous presence in the formulae for various exit probabilities.

In particular, if the Brownian motion with drift $\brmd{m}$ starts at $s$, the probability that it exits the interval $[s-\Delta_l,s+\Delta_r]$ at the right (left) end is, respectively,
      \be\label{eq:exit}
      p(\Delta_l,\Delta_r)=\frac{1-\Em(\Delta_l)}{1-\Em(\Delta_l+\Delta_r)}; q(\Delta_l,\Delta_r)=\frac{\Em(\Delta_l)-\Em(\Delta_l+\Delta_r)}{1-\Em(\Delta_l+\Delta_r)}.
      \ee

\subsubsection{Automata and Sums over Paths}
      Consider a general finite state automaton, with edges marked by weights: certain monomials in some collection of formal variables $z_1, z_2,\ldots, z_m$ with real coefficients. 
      We will denote the weight corresponding to the edge $(\beta,\beta')$ as $w(\beta,\beta')$. For a finite path $\pth$ in the automaton, we define the weight $w(\pth)$ of the path to be the product of all the weights of the transitions.

      For a pair of states $\beta,\beta'$ in the automaton, consider the formal power series (with nonnegative real coefficients) formed by a sum over all paths going from $\beta$ to $\beta'$:
      \be\label{eq:sum_paths}
      F_{\beta,\beta'}(z)=\sum_{\arpath{\beta}{{\pi}}{{\beta}'}} w(\pth),
      \ee

      The usual rules apply, for example, the  the most important being the {\em last edge decomposition},
      \be\label{eq:rule}
      F_{\beta,\beta''}=\delta_{\beta,\beta''}+\sum_{\beta\to\beta'} F_{\beta,\beta'}w(\beta',\beta'').
      \ee
This leads, as usual, to rational generating functions for the enumerations of weights over trajectories in the automaton.

\subsubsection{Distributions of the Windings around Intervals}
In the situation where the edge weights are transition probabilities multiplied by monomials in formal variables, the sums over paths \eqref{eq:sum_paths} become partition functions of the various descriptors of the trajectories. If the trajectory spends only finite time in the span of the automaton, then the sum \eqref{eq:sum_paths} has well-defined terms, whose nonnegative coefficients sum up to the unity (and converging in the unit polydisk $\{|z_l|\leq 1, l=1,\ldots,m\}$), and the standard computations render each of $F_{\beta,\beta'}(z)$ a rational function of the variables $z_,\ldots, z_m$.
      
As an example, if one associates to the automaton shown on Fig. \ref{fig:auto_a} transition probabilities
\[
p_{\delta,\beta}=\Em(\Delta); p_{\delta,\omega}=1-\Em(\Delta), \mbox{ where } \Delta=d-b; \mbox{ all others }=1,
\]
and sets the weight on $(\delta,\beta)$ equal to $p_{\delta,\beta}z$, then
\[
F_{\alpha,\omega}=\sum_k p_k z^k,
\]
where $p_k$ is the probability that a trajectory of the Brownian motion with drift has $k$ finite bars spanning $[b,d]$. Applying \eqref{eq:rule} we obtain
\[
F_{\alpha,\alpha}=1; F_{\alpha,\beta}=F_{\alpha,\alpha}+z \Em(\Delta)F_{\alpha,\delta}; F_{\alpha,\delta}=F_{\alpha,\beta}; F_{\alpha,\omega}= (1-\Em(\Delta))F_{\alpha,\delta},
\]
resolving which leads, unsurprisingly, to 
\[
F_{\alpha,\omega}=\frac{1-\Em(\Delta)}{1-z\Em(\Delta)}.
\]

\subsection{Higher Windings}\label{ss:higher}
Understanding joint windings around several intervals requires construction of somewhat more complicated automata. 

One useful construction is that of the {\em product} of two automata.

\subsubsection{Product of Automata}
Consider a pair of automata of the kind described in Section \ref{sss:automata}: two collections of intervals, $\intervals_1, \intervals_2$, and the corresponding switch maps $s_1, s_2$ defined on the endpoints $\boundaries_1,\boundaries_2$. Assume, for simplicity, that the heights of the finite ends of the intervals in different collections are all disjoint.

Then one can form a new automaton, {\em the product of the original ones}. The intervals of the product automaton would be the formed by the pairs $U_{\inl_1,\inl_2}=U_{\inl_1}\cap U_{\inl_2}$ for which this overlap is nonempty. Note that heights of different intervals of this kind can coincide: those intervals are enumerated not by their heights, but by their indices $(\inl_1,\inl_2)$.

Each end of the interval $U_{\inl_1,\inl_2}$ comes unambiguously from one if the collections, $\intervals_1$ or $\intervals_2$. This would allow to define the switch map on the product of the original automata. If, say, the upper end of the interval $(\inl_1,\inl_2)^+$ stems from the first automaton (that is the upper height of $\inl_1$ is less than that of $\inl_2$), then the switch map $s$ of the sends $(\inl_1,\inl_2)^+$ to $(\tau_1(\inl_1^+),\inl_2)$, and so on.

The Fig. \ref{fig:higherw} illustrates the products of two elementary automata for differently overlapping intervals.

We remark that in the case when the intervals are not nested, one of the branches is a {\em "garden of Eden"}, - a trajectory starting outside it will never enter it. Hence one can ignore it; this reduces the effective number of states to six.

\subsubsection{Windings around a Pair of Intervals}
     Let $I_1=[b_1,d_1], I_2=[b_2,d_2]$ are two intervals. We will be looking into the joint distribution of the numbers $b_1, b_2$ of bars spanning these intervals.

\begin{figure}[htbp]
        \captionsetup{font={small, it},width=.8\textwidth}
        \begin{center}
          \hfill   \includegraphics[height=1.85in]{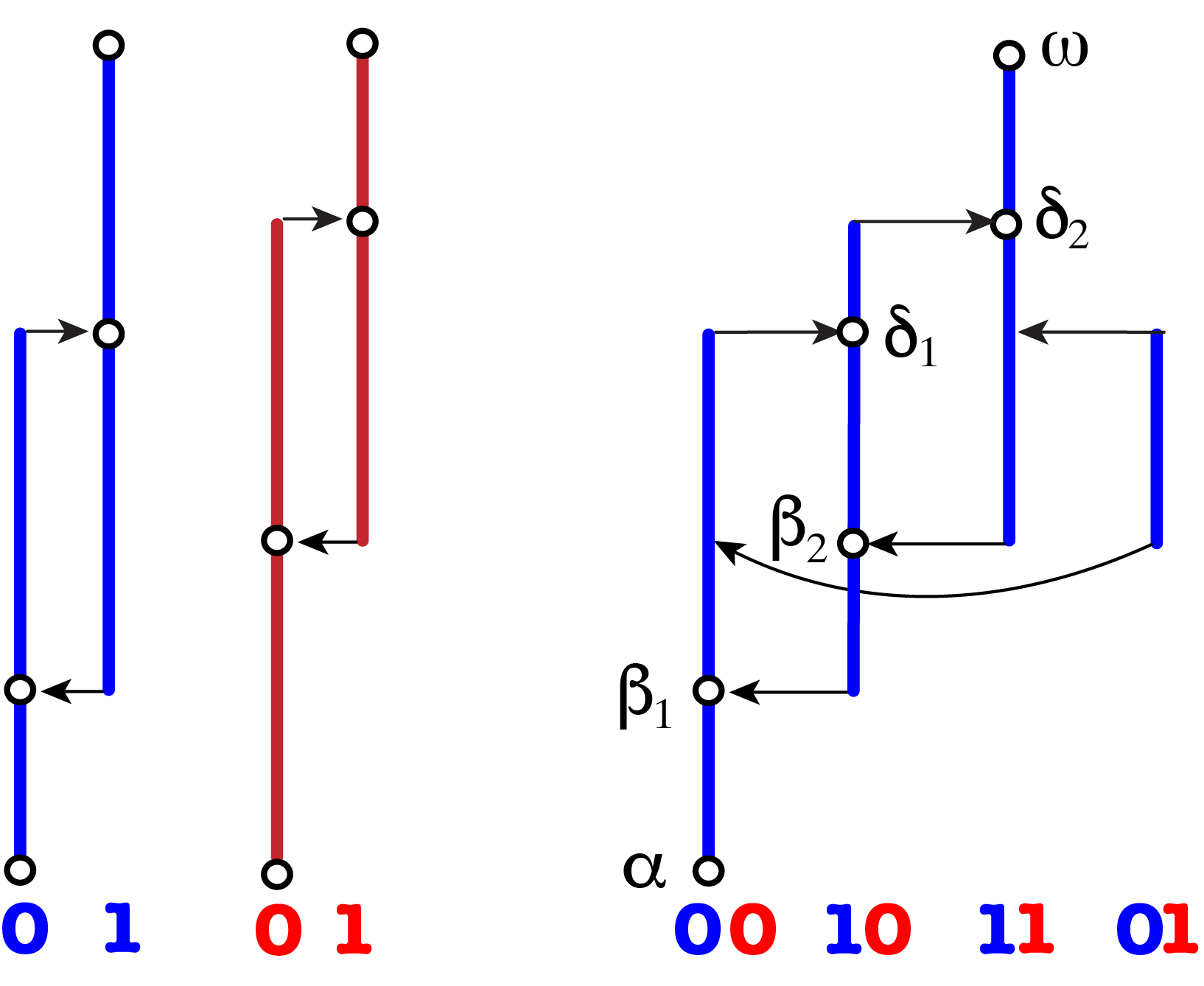}\hfill   
          \includegraphics[height=1.7in]{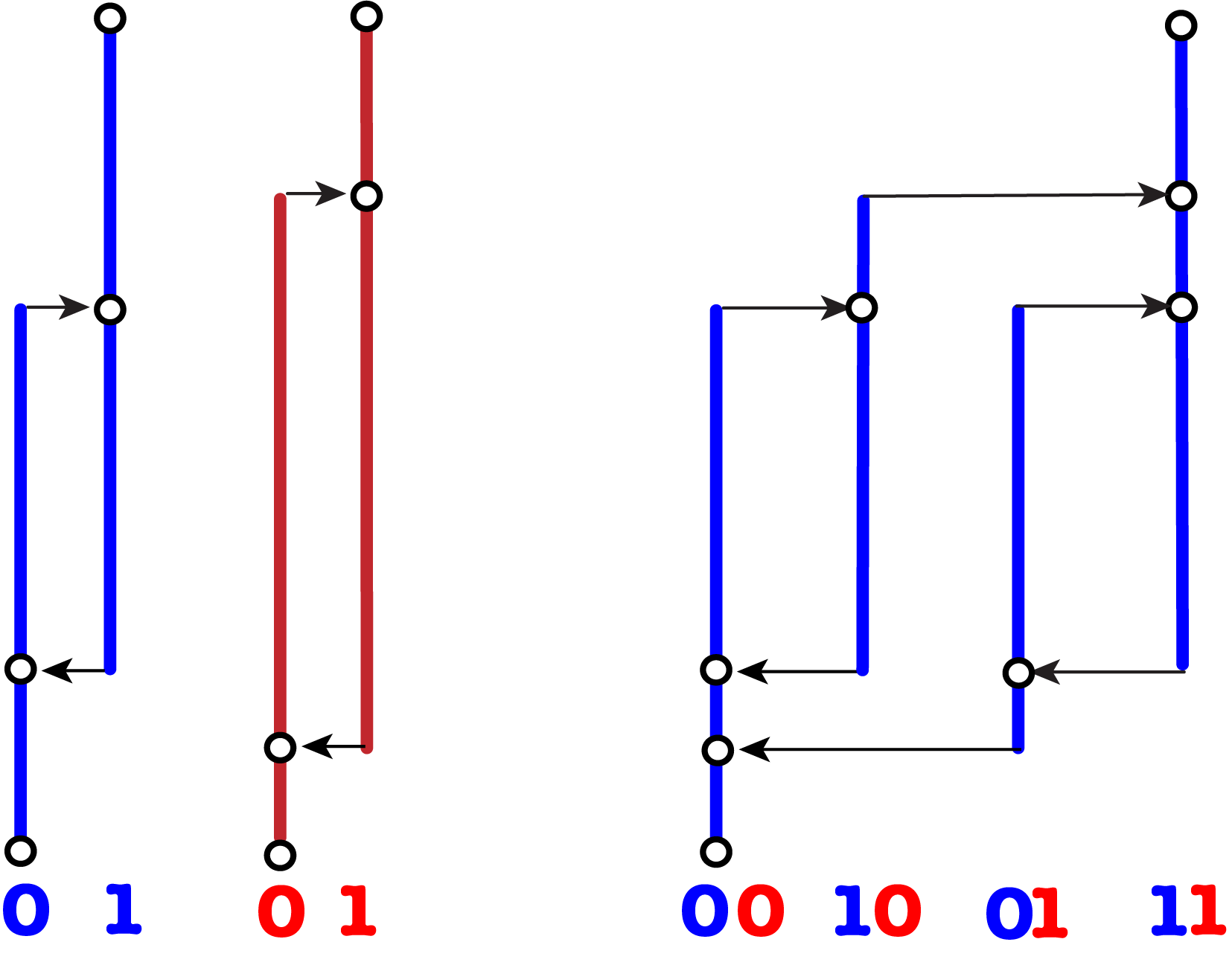}\hfill\hfill 
          \caption{Windings around a pair of intervals. We do not indicate the states corresponding to the exits from the interval {\tt 01} on the left, as they do not contribute to the final result.}
          \label{fig:higherw}
        \end{center}
      \end{figure}

Assume that the intervals are not nested, so that $b_1<b_2; d_1<d_2$. In this case, the relevant automaton is shown on the left display of Fig. \ref{fig:higherw}. Any continuous function that starts below $b_1$ and ends above $d_2$ generates a symbolic trajectory. The following is immediate:

      \begin{lem}
        The total number of windings around the interval $I_1$ is counted by the number of transitions into the state $\beta_1$,  (i.e. either $\delta_1\to\beta_1$ or $\beta_2\to\beta_1$), while the number of windings around $I_2$ is given by the number of transitions $\delta_2\to\beta_2$.
      \end{lem}
\begin{proof}
Indeed, lifts of the trajectory to either of the factors coincides with the projection to that factor of the lift of the trajectory into the product, and the transitions in either lift happen at the same positions.
\end{proof}

Let us introduce the generating functions for the numbers of windings around the intervals $I_1,I_2$ for the paths starting at $\alpha$ and ending at a state $\nu$:
      \be\label{eq:gen_fun_a}
      H_\nu(x,y):=F_{\alpha,\nu}(x,y)=\sum_{\arpath{\alpha}{\pi}{\nu}} w(\pth)=\sum_{\arpath{\alpha}{\pth}{\nu}}x^{k(\pth)} y^{l(\pth)}\prob_{\pth}:
      \ee
here $k(\pth),l(\pth)$ are numbers of windings of the path $\pth$ around the intervals $I_1$, $I_2$, respectively, $\prob(\pth)$ is the probability of the realization of the symbolic path $\pth$, and the summation is over all paths starting in $\alpha$ and ending in $\nu$.

We will conduct the computations only for the non-nested intervals, depicted on the left display of Fig. \ref{fig:higherw}. The case of nested intervals, as well as more correlation functions of higher orders will be addressed elsewhere.

Using the last edge decomposition \eqref{eq:rule} for the chain on the Figure \ref{fig:higherw} we arrive at
      \be
      \begin{array}{rl}
        H_\alpha&=1\\
        H_{\beta_1}&=H_\alpha+x(q_{\delta_1}H_{\delta_1}+q_{\beta_2}H_{\beta_2})\\
        H_{\beta_2}&=y q_{\delta_2}H_{\delta_2}\\
        H_{\delta_1}&=H_{\beta_1}\\
        H_{\delta_2}&=p_{\delta_1}H_{\delta_1}+p_{\beta_2}H_{\beta_2}\\
        H_\omega&=p_{\delta_2} H_{\delta_2}.\\
      \end{array}
      \ee
Here we use the shortcuts $p_\nu, q_\nu$ to denote the probabilities to exit from the top or bottom end, respectively, from the interval containing the state $\nu$, conditioned on the Brownian trajectory starting at $\nu$.

Now we are ready to compute the joint distributions of the numbers of windings around the intervals $I_1,I_2$.

      The formal power series \eqref{eq:gen_fun_a} is, clearly, the generating function for the numbers of windings around the intervals $I_1, I_2$.
\begin{prp}
The generating function $H:=H_\omega$ is given by 

      \be\label{eq:H}
      H(x,y)=-\frac{(\Eb_1-\Ed_1) (\Eb_2 - \Ed_2)}
      {\Ed_1(\Eb_1-\Ed_2) + \Eb_1(\Ed_2-\Ed_1) x +
      \Ed_1( \Eb_2- \Eb_1) y +\Eb_1(\Ed_1-\Eb_2) x y}.
      \ee
      where we denote $\Eb_k=\Em(-b_k);\Ed_k=\Em(-b_k), k=1,2$.
\end{prp}
\begin{proof}
The work is shown in the Appendix B.
\end{proof}
      \subsection{Correlation density}
      Using equation \eqref{eq:H}, we can easily find the covariance of the numbers $\wi_1,\wi_2$ of windings around the intervals $I_1,I_2$: indeed,
      \be\label{eq:cov}
      \cov(I_1,I_2)=\left.\ex \wi_1\wi_2 -\ex \wi_1 \ex \wi_2=\frac{\partial^2 \log H}{\partial x\partial y}\right\vert_{x=y=1}.
      \ee

      Performing the calculations results in
      \be\label{eq:covkl}
      \cov(I_1,I_2)=\frac{\exp(-2m(d_2-b_1))}
      {(\exp(-2m(d_1-b_1))-1)(\exp(-2m(d_2-b_2))-1)}.
      \ee

Recall that any point process is characterized by its {\em factorial moments}, which in our tame setting can be characterized by their densities, i.e. {\em correlation functions}.

\begin{prp}
The second moment density function of $\ph_0$ for the non-nested intervals (i.e. for $b_1<b_2, d_1<d_2$) is given by
      \be
      \mph_2(b_1,d_1,b_2,d_2)=\frac{\partial^4 \cov(b_1,d_1,b_2,d_2)}
      {\partial b_1\partial d_1\partial b_2\partial d_2}.
      \ee

\end{prp}

\begin{proof}
The covariance \eqref{eq:covkl} can be interpreted as a function of the endpoints $(b_1,d_1,b_2,d_2)$ as follows:
\be\label{eq:bi}
\cov(b_1,d_1,b_2,d_2)=\\ \ex\left(\int {\bm \phi}_{b_1,d_1}d\ph_0 \int{\bm \phi}_{b_2,d_2}d\ph_0\right)-
\ex \left(\int {\bm \phi}_{b_1,d_1}d\ph_0\right) \ex\left(\int{\bm \phi}_{b_2,d_2}d\ph_0\right),
\ee
where ${\bm \phi}_{b,d}$ is the indicator function of the NW coordinate quadrant with the apex at $(b,d)$. 

The bilinearity of \eqref{eq:bi} with respect to ${\bm{\phi}}$'s implies that the covariance of the $\ph_0$-contents of the small coordinate aligned squares centered at $(b_1,d_1)$ and $(b_2,d_2)$ can be obtained by taking first order differences with respect to four variables $(b_1,d_1,b_2,d_2)$ of \eqref{eq:covkl}. 

Dividing by the areas of these squares and passing to the limit, we obtain \eqref{eq:bi}
\end{proof}

Lengthy but elementary computations yield
\begin{thm}\label{thm:g2}
For non-nested intervals $b_1<b_2; d_1<d_2$, the 2-point correlation function for the point process $\ph_0$ of 0-dimensional persistence for the standard Brownian motion with drift $m$ is given by
      \be\label{eq:covdef}
      \mph_2(b_1,d_1;b_2,d_2)=\frac
      {64 m^4 \exp(2 b_1 + b_2 + 2 d_1 + d_2)}
      {(\exp(2m b_1) - \exp(2m d_1))^3 (\exp(2m b_2) - \exp(2m d_2))^3}.
      \ee
\end{thm}

The expression \eqref{eq:covdef} can be rendered as
\[\mph_2(b_1,d_1;b_2,d_2)=
\frac
      {64 m^4 \Em(\Delta_1)\Em(\Delta_2)\Em(\Delta)}
      {(1-\Em(\Delta_1))^3 (1-\Em(\Delta_2))^3}=
      \mph(\Delta_1)\mph(\Delta_2)\frac{4\Em(\Delta)}{(1+\Em(\Delta_1))(1+\Em(\Delta_2))};
\]
here $\mph$ is the intensity density of the point process $\ph_0$ (see the Proposition \ref{prp:intensity}), $\Delta_1, \Delta_2$ are the lengths of the intervals $I_1, I_2$ and $\Delta=d_2-b_1$ the span of the union of the intervals.

\begin{rem} The Theorem \ref{thm:g2} implies, in particular, that the {\em pair correlation function} is positive, at least for the non-nested intervals. 

Another implication is the exponential decay of correlation function. By standard arguments, this implies that the expected number and the variance of the number of bars in the region 
\[
Q_c(L):=\{0<b<d<L, d-b>c>0\}
\]
grow asymptotically linearly with $L$, and that the number itself is asymptotically normal, as $L\to\infty$.
\end{rem}

\section{Concluding remarks}
\begin{itemize}
\item As we mentioned above, it would be of interest to understand the overall structure of the point process, in particular to derive its higher factorial moments. The techniques of the products of elementary automata makes the problem computable.

\item Similarly, the question of the {\em joint distribution} of the numbers of bars of either type, \fr and \maa,  is of interest. 

\item We focused in this note on the special class of Brownian motions. One might expect some universality near the diagonal. For example, the limits of the small drift can be obtained using the rescaling $t\mapsto \lambda t; h\mapsto \lambda^{1/2}h$ and the invariances described in Section \ref{sec:invar}.
\item An alternative approach to description of the persistence diagrams of trajectories of rather general $\Real$-valued Markov random processes was addressed by J. Picard in \cite{picard_tree_2008}: he generated a sequence of stopping times marking emergence of all bars of length above some threshold, and used its properties to derive quite general results on the resulting point process. This line of research was also pursued by D. Perez, see e.g. \cite{perez}.

\item As we mention in Appendix A, the behavior of the stack sizes in the algorithm for finding the barcodes of a time series seems nontrivial and interesting, with random {\em updown permutation} \cite{arnold_calculus_1992} being the natural model for the algorithm input.

\item Generalizing to higher-dimensional domains seems difficult but not entirely out of reach, at least for some classes of Gaussian processes, using the techniques on the bar decompositions of the products of random merge trees. This line of research will be pursued elsewhere.

\end{itemize}

\subsection*{Acknowledgements}
This research started from a conversation with Robert Adler, Jonathan Taylor and Shmuel Weinberger. 

      \section{Appendices}

      \subsection*{A. Finding the Bars}\label{sec:algorithm}
      Finding the bars of a time series can be done quite efficiently.

      Assume that a function is given as a time series, i.e. a list of its values $f(k), k=1,\ldots,N$ at consecutive points $t_1<t_2<\ldots <t_N$ between which the function is assumed to be monotonic (say, linearly interpolating). Further we will assume that the function is generic.

      We augment the series by setting $f(0)$ and $f(N+1)$ to respectively global minimum and maximum values: $f(0)< f(k)<f(N+1), k=1,\ldots,N$.

      The output is a list with entries having the structure
      $(d,t_d;b,t_b)$, where $(d<b)$ is a bar, and $(t_d, t_b)$ are the
      locations of the corresponding (local) maxima and minima.

      The algorithm maintains at all times two stacks, one with local minima, one with local maxima, ordered by critical value within each stack. Removals from the stacks happen in pairs and produce bars.

      \begin{algorithm}{\small
        \caption{One pass bar algorithm for time series.}
        \begin{algorithmic}
          \State $\Top,\Bot\gets \mathtt{empty\ stacks}$ \Comment{Initializing}
          \State $\Dir\gets+1$ \Comment{Initially, the function increases}
          \For {$t=1,\ldots,N$}
          \If {$(f(t)-f(t-1))\times \Dir <0$} \Comment{Direction changes, so either}
          \If {$\Dir=+1$}
          \State {$\Top.\Push((t-1,f(t-1))$} \Comment{add the local maximum}
          \Else
          \State {$\Bot.\Push((t-1,f(t-1))$} \Comment{...or minimum to respective stack,}
          \EndIf
          \State{$\Dir=-\Dir$} \Comment{and record change of the direction.}
          \Else
          \If {$(\Dir=+1 \&\& f(t)>\Top(0))||(\Dir=-1 \&\& f(t)<\Bot(0))$}
          \State {$\mathtt{Output}(\Top.\Pop,\Bot.\Pop)$}\Comment{If hitting a wall, output a bar.}
          \EndIf
          \EndIf
          \EndFor
        \end{algorithmic}}
      \end{algorithm}


While the time execution of the algorithm is clearly linear in the length of the time series, an interesting question arises on the depth of stacks (memory) required for its execution. In the worst case (when the local maxima are descending and the local minima ascending) the depth required is linear as well. What about the average case?

We remark that one can adapt this algorithm to register arbitrarily complicated snake-like patterns in time series.

\subsection*{B. Merge Tree Monoid}
Consider two path-connected spaces with continuous functions on them: $f:X\to \Real, g:Y\to \Real$. We define the Thom-Sebastiani sum of the functions $f$ and $g$ as the function $f\oplus g:X\times Y\to \Real$ taking $f\oplus g:(x,y)\mapsto f(x)+g(y)$ (motivated by \cite{thom_sebastiani}).

One can prove
\begin{prp}
The merge tree $\mtree_{f\oplus g}$ of the Thom-Sabatini sum depends only on the merge trees $\mtree_f, \mtree_g$.
\end{prp}

In particular, one can consider $\oplus$ as an operation on the merge trees (with roots placed at the infinite height).

\begin{dfn}
The operation $\oplus$ turns the collection of merge trees bounded from below into a {\em commutative monoid}, with the unity given by the merge tree $\mtree_\zero$ whose sole leaf is at height $0$ (the merge tree of any positive definite quadratic function, analogously to Thom-Sabatini setup).
\end{dfn}

We note that the merge trees corresponding to $x^2+r$, i.e. tree with a single leaf at the height $r$ form an embedding of $\Real$ into the merge trees monoid, and addition with them realizes the shift action of $\Real$ on the monoid.

  To recover $\mtree_{f\oplus g}$ from the merge trees $\mtree_f, \mtree_g$, up to chiralities, it is enough to identify all bars (i.e. stems resulting from the recursive pruning of the tree, see Section \ref{sss:mttoph}) and their attachments. The collection of bars of $\mtree_{f\oplus g}$ is easy to recover: they are in one-to-one correspondence to pairs of bars $[b_f,d_f]$, $[b_g,d_g]$ of $f$ and $g$, and have birth at $b_f+b_g$, and have the lifespan equal to the smaller of the lifespans of $d_f-b_f, d_g-b_g$.

\begin{figure}[htbp]
  \captionsetup{font={small, it},width=.8\textwidth}
        \begin{center}
          \includegraphics[height=2.6in]{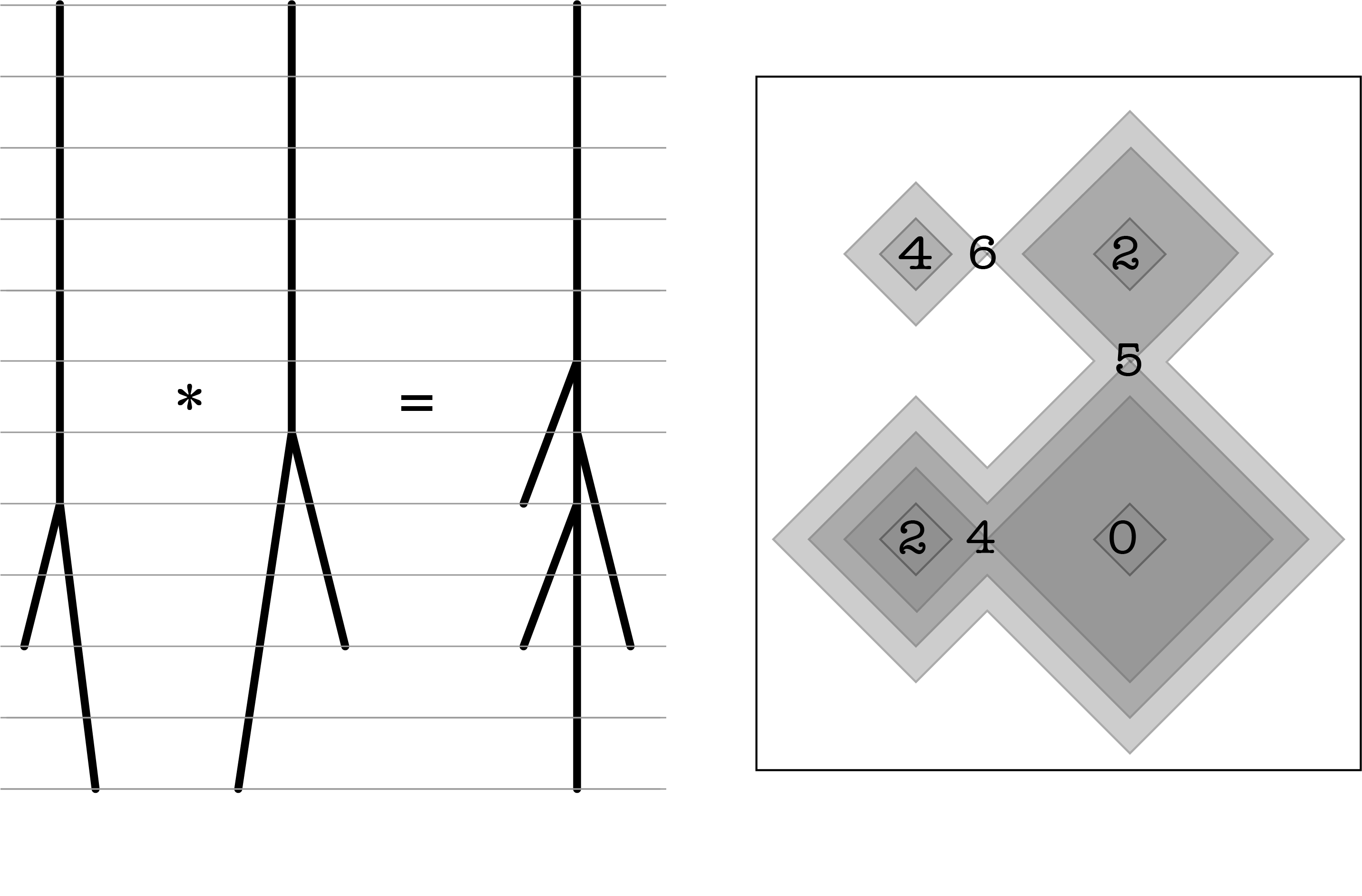}          
          \caption{Product of merge trees.}
          \label{fig:merge}
        \end{center}
      \end{figure}

Consider two merge trees shown on Fig. \ref{fig:merge} together with their product; the right display shows some of the level sets of $f(x)+g(y)$, where $f,g$ are the univariate functions whose merge trees correspond to the factors on the left.

The rules of the attachments, details of the construction, as well as applications to the study of persistence diagrams of (tilted) Brownian sheets will be presented in a forthcoming publication \cite{PLe}.

Another application of merge tree monoid structure could be the model of {\em corner percolation}, see \cite{pete_corner_2008}. There one studies the structure of the level sets of functions $f(x)-g(y)$ where $f, g$ are (univariate) random functions representing the trajectories of random excursions of equal height.

      \subsection*{Appendix C: Generating functions for windings around two intervals}
      Solving for $F_\omega$ is routine, the result is
      \be\label{eq:gencor}
      H(x,y)= p p_1/(1 - x q_1 - y q p_2 + x y q (q_1 p_2 - p_1 q_2))
      \ee
      where we use shorthands (relying on \eqref{eq:exit})
      \begin{eqnarray*}
        p=p( d_2- b_2,\infty);\hfill q=q( d_2- b_2,\infty);\\
        p_1=p( d_1- b_1, d_2- d_1); q_1=q( d_1- b_1, d_2- d_1);\\
        p_2=p( b_2- b_1, d_2- b_2); q_2=q( b_2- b_1, d_2- b_2).
      \end{eqnarray*}
      Using the expressions \eqref{eq:exit} for $p$'s and $q$'s we find that $H(x,y)$ is given by

      \begin{equation*}
        -\frac{(\Em( b_1)-\Em( d_1)) (\Em( b_2) - \Em( d_2))}
        {\Em( d_1)(\Em( b_1)-\Em( d_2)) + \Em( b_1)(\Em( d_2)-\Em( d_1)) x +
        \Em( d_1)( \Em( b_2)- \Em( b_1)) y +\Em( b_1)(\Em( d_1)-\Em( b_2)) x y},
      \end{equation*}
which simplifies to \eqref{eq:H}.

\section{Declarations}
 
\subsection{Ethical Approval}
Not applicable.
 
\subsection{Funding}
Partially supported by the NSF via
  grant DMS-1622370
 
\subsection{Availability of data and materials}
Not applicable

      \bibliographystyle{alpha}
      \bibliography{r}

      \end{document}